# $\vec{C}$ Sequential Optimization Numbers Group

Zile Hui, 51174500096@stu.ecnu.edu.cn

ABSTRACT. We define **C** sequential optimization numbers, where **C** is a *k+1*-tuple vector. We prove that the unsigned Stirling numbers of first kind are (0,1) sequential optimization numbers. Many achievements of the Stirling numbers of first kind can be transformed into the properties of **C** sequential optimization numbers. We give some examples such as the recurrence formula and an instance of **C** sequential optimization numbers. We also extend some properties such as an upper bounder of them.

## 1 INTRODUCTION

Stirling numbers were introduced by the Scottish mathematician James Stirling in his famous treatise and were subsequently rediscovered in various forms by numerous authors [Bla16][Sti30]. Stirling numbers of the first kind, denoted by $s(n, m)$, were studied in a large number of works and have various properties [Bla16]. Recurrence formula is one of the basic ways to define the unsigned Stirling numbers of first kind. There is no known closed-form expression so far for the unsigned Stirling numbers of first kind and asymptotic formulas have been studied by numerous authors [Wil93][Tem93][AD17]. **C** Sequential optimization numbers are defined based on optimization set. We give some properties of **C** sequential optimization numbers and some of them are fundamental properties of Stirling numbers when **C** is (0,1). We mainly refer to the definitions and proofs in one article [Hui22].

## 2 SEQUENTIAL OPTIMIZATION NUMBERS

We introduce the optimization set in Definition 2.1 and define the **C** sequential optimization numbers in Definition 2.2. Then, we give an expression which is not closed-form, a recurrence formula, some properties and applications of **C** sequential optimization numbers. At the end of this chapter, we give an upper bound of sequential optimization numbers and some results derived from this upper bound.

**Definition 2.1.** *Let* $\boldsymbol{x}(x_1, x_2, \dots, x_k)$, $\boldsymbol{y}(y_1, y_2, \dots, y_k)$ *be k-dimensional vectors and* $\boldsymbol{R}(R_1, R_2, \dots, R_k)$ *be k-dimensional relation vector. If for all* $i = 1,2, \dots, k$, $(x_i, y_i) \in \boldsymbol{R_i}$, *then* $\boldsymbol{x}$ *is said to be related to* $\boldsymbol{y}$ *by* $\boldsymbol{R}$, *denoted by* $(\boldsymbol{x}, \boldsymbol{y}) \in \boldsymbol{R}$. *Let U be a set of vectors,* $\boldsymbol{R}$ *be a relation vector and* $A \subseteq U$, *if for* $\forall \boldsymbol{u} \in U$, $\exists \boldsymbol{a} \in A$ *imply* $(\boldsymbol{a}, \boldsymbol{u}) \in \boldsymbol{R}$ *or* $\boldsymbol{a} = \boldsymbol{u}$, *then A is said to be a majorization set of U by* $\boldsymbol{R}$. *If B is a majorization set of U by* $\boldsymbol{R}$ *and* $|B|$ *is a minimum of all the cardinalities of majorization sets, then B is said to be an optimization set of U by* $\boldsymbol{R}$, *denoted by* $B \overset{o\boldsymbol{R}}{\subseteq} U$, *and* $|B|$ *is said to be the weight of U by* $\boldsymbol{R}$, *denoted by* $W_{o\boldsymbol{R}}(U)$.

We can find similar selection strategies of optimization set in many existing studies [Mou79][Xue00][ARY21][Hui22]. It can be proved that this strategy can include all possible optimal terms and exclude impossible optimal terms in the iteration of the algorithm [Hui22].

**Definition 2.2.** *Let* $\boldsymbol{a_1}(1, a_{11}, a_{12}, \dots a_{1k})$, $\boldsymbol{a_2}(2, a_{21}, a_{22}, \dots a_{2k})$, …, $\boldsymbol{a_n}(n, a_{n1}, a_{n2}, \dots a_{nk})$ *be n k+1-dimensional vectors,* $\boldsymbol{R}(R_1, R_2, \dots, R_k) = (<, <, \dots, <)$ *be k-dimensional relation vector,* $\boldsymbol{C} = (c_0, c_1, c_2, \dots, c_k)$, $U = \{\boldsymbol{a_1}, \boldsymbol{a_2}, \dots, \boldsymbol{a_n}\}$ *and* $S \subseteq U$. *For all* $j = 0,1, \dots, k$, $c_j \in \{0,1\}$. *For all* $i = 1,2, \dots, n$ *and* $j = 1,2, \dots, k$, $\boldsymbol{b_{ij}} = (i, a_{ij})$, $U_j = \{\boldsymbol{b_{1j}}, \boldsymbol{b_{2j}}, \dots, \boldsymbol{b_{nj}}\}$, $a_{1j}, a_{2j}, \dots a_{nj}$ *are* $1,2, \dots, n$ ,

*respectively.* $S_j$ *is an optimization set of* $U_j$ *by* $(<, R_j)$*. For a certain* $i$ *and all* $j = 1,2, \ldots, k$*, the total number of* $\boldsymbol{b_{ij}} \in S_j$ *is* $l$*, if* $c_l = 1$*, then* $\boldsymbol{a_i} \in S$*, otherwise,* $\boldsymbol{a_i} \notin S$*,* $S$ *is said to be an* $\boldsymbol{C}$ *sequential optimization set of* $U$ *by* $\boldsymbol{R}$*, denoted by* $S \overset{\boldsymbol{OR}, \boldsymbol{C}}{\subseteq} U$*.* $|S|$ *is said to be sequential optimization weight of* $U$ *by* $\boldsymbol{R}$*, denoted by* $W_{\boldsymbol{OR},\boldsymbol{C}}(U)$*. The numbers of ways that* $W_{\boldsymbol{OR},\boldsymbol{C}}(U) = m$ *are said to be* $\boldsymbol{C}$ *sequential optimization numbers, denoted by* $O_{\boldsymbol{C}}(n, m)$*.*

In Definition 2.2, $\boldsymbol{C} = (c_0, c_1, c_2, \ldots, c_k)$. We define $\boldsymbol{B}_j = (b_{j,0}, b_{j,1}, b_{j,2}, \ldots, b_{j,k})$, $\boldsymbol{C}' = (1,1, \ldots, 1) - \boldsymbol{C}$ and $F_j(\boldsymbol{C}) = \boldsymbol{B}_j \cdot \boldsymbol{C}^T$, where $b_{j,p} = C_k^p \frac{1}{(j-1)^p}$, $p = 0,1, \ldots, k$ and $j = 2,3, \ldots, n$. In this paper, we define $\prod_{j=m+1}^{m} F_j(\boldsymbol{C}) = 1$ and $\sum_{p=m+1}^{m} f(p) = 0$ to make the expressions concise, where $f(p)$ is a function of $p$. We give some properties about $F_j(\boldsymbol{C})$.

*(a). If* $2 \le j_1 \le j_2 \le n$*,*

$$F_{j_2}(\boldsymbol{C}) \le F_{j_1}(\boldsymbol{C}) \tag{1}$$

**Proof.** If $2 \le j_1 \le j_2 \le n$,

$$\frac{b_{j_2,p}}{b_{j_1,p}} = \frac{(j_1 - 1)^p}{(j_2 - 1)^p} \le 1$$

$$b_{j_2,p} \le b_{j_1,p}$$

$$F_{j_2}(\boldsymbol{C}) \le F_{j_1}(\boldsymbol{C})$$

(b).

$$\prod_{j=2}^{n} F_j(\boldsymbol{C}) \le n^k \tag{2}$$

**Proof.**

$$\begin{aligned}\prod_{j=2}^{n} F_j(\boldsymbol{C}) &\le \prod_{j=2}^{n} F_j([1,1, \ldots, 1]^T) \\ &= \prod_{j=2}^{n} \sum_{p=0}^{k} b_{j,p} \\ &= \prod_{j=2}^{n} \sum_{p=0}^{k} C_k^p \frac{1}{(j-1)^p} \\ &= \prod_{j=2}^{n} \left(1 + \frac{1}{j-1}\right)^k \\ &= n^k\end{aligned}$$

*(c). If* $c_0 = 1$*,*

$$F_j(\boldsymbol{C}) \ge 1 \tag{3}$$

**Proof.**

$$F_j(\boldsymbol{C}) = \boldsymbol{B} \cdot \boldsymbol{C}^T \ge c_0 b_{j,0} = 1$$

**Theorem 2.1.** *For all* $k, n \in N^+$ *and* $m \in N$*,*

$$O_{\boldsymbol{C}}(n, m + c_k - 1) = \begin{cases} (n-1)!^k \sum \left( \prod_{i=1}^{m-1} F_{j_i}(\boldsymbol{C}) \prod_{i=1}^{n-m} F_{j'_i}(\boldsymbol{C}') \right), & 1 \le m \le n \\ 0, & otherwise \end{cases}$$

*where* $j_1, j_2, \cdots, j_{m-1}$ *are all combinations consisting of m-1 elements in set* $\{2,3, \ldots, n\}$ *and* $\{j'_1, j'_2, \cdots, j'_{n-m}\} = \{2,3, \ldots, n\} - \{j_1, j_2, \cdots, j_{m-1}\}$*.*

**Proof of Theorem 2.1.** In the scenario of Definition 2.2, $\boldsymbol{C} = (c_0, c_1, c_2, \dots, c_k)$ and $S \overset{\boldsymbol{OR},\boldsymbol{C}}{\subseteq} U$. $g(S)$ denotes the numbers of ways where $S$ is a set consisting of specific elements.

First, we discuss the case where $\boldsymbol{C} = (0,1)$. We can get $S$ is optimization set.

For $m = 1$, the number of positions of $\boldsymbol{a_2} \sim \boldsymbol{a_n}$ is $(n-1)!$ as shown in Figure 1.*a*. So $O_{(0,1)}(n, 1) = (n-1)!$.

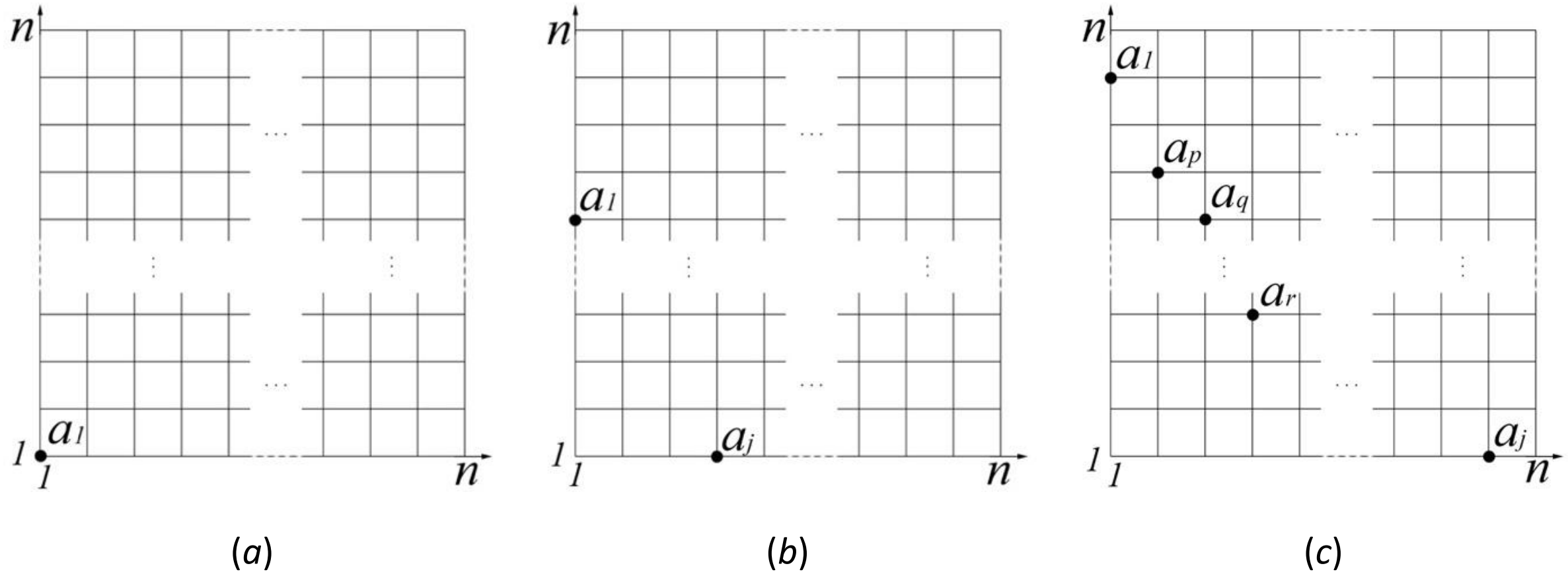

Figure 1: Three ways to optimization set.

For $m = 2,3,\dots,n$, $|U| = n$ and $S_1 = \{\boldsymbol{a_1}, \boldsymbol{a_p}, \boldsymbol{a_q}, \dots, \boldsymbol{a_r}, \boldsymbol{a_j}\}$, we use mathematical induction to prove that $g(S_1) = \frac{(n-1)!}{(p-1)(q-1)\cdot\ldots\cdot(r-1)(j-1)}$.

Basis step: we prove $g(S_2) = \frac{(n-1)!}{j-1}$, where $S_2 = \{\boldsymbol{a_1}, \boldsymbol{a_j}\}$. As shown in Figure 1.*b*, $a_{11} = i$, $a_{j1} = 1$, where $i = 2,3,\dots,n$, $j = 2,3,\dots,n$ and $i + j \le n + 2$. First, place row $2\sim(i-1)$ and the number of positions is $\frac{(n-j)!}{(n-i-j+2)!}$. Second, place row $(i+1)\sim n$ and the number of positions is $(n-i)!$. Total number is $\frac{(n-i)!(n-j)!}{(n-i-j+2)!}$.

$$g(S_2) = \sum_{i=2}^{n+2-j} \frac{(n-i)!\,(n-j)!}{(n-i-j+2)!}$$

$$= \frac{(n-j)!}{j-1} \sum_{i=2}^{n+2-j} \frac{(n-i)!\,[(n-i+1)-(n-i-j+2)]}{(n-i-j+2)!}$$

$$= \frac{(n-j)!}{j-1} \sum_{i=2}^{n+2-j} [\frac{(n-i+1)!}{(n-i-j+2)!} - \frac{(n-i)!\,(n-i-j+2)}{(n-i-j+2)!}]$$

In the brackets, the right side of the minus is equal to 0 when $i = n + 2 - j$. The right side of the minus with $i = t$ is equal to the left side of the minus with $i = t + 1$, where $t = 2,3,\dots,n+1-j$. So,

$$g(S_2) = \frac{(n-j)!}{j-1} \cdot \frac{(n-1)!}{(n-j)!}$$

$$= \frac{(n-1)!}{j-1}$$

Inductive step: assume $g(S_3)=\frac{(n-1)!}{(p-1)(q-1)\cdot\ldots\cdot(r-1)}$, where $S_3=\{\boldsymbol{a_1},\boldsymbol{a_p},\boldsymbol{a_q},\ldots,\boldsymbol{a_r}\}$. Then, we show $g(S_1)=\frac{(n-1)!}{(p-1)(q-1)\cdot\ldots\cdot(r-1)(j-1)}$ where $S_1=\{\boldsymbol{a_1},\boldsymbol{a_p},\boldsymbol{a_q},\ldots,\boldsymbol{a_r},\boldsymbol{a_j}\}$ and $1<p<q<r<j\le n$. In Figure 1.*c*, $a_{r1}=i$, $a_{j1}=1$, where $i=2,3,\ldots,n-1$, $j=3,4,\ldots,n$ and $i+j\le n+2$. First, place row $2\sim(i-1)$ and the number of positions is $\frac{(n-j)!}{(n-i-j+2)!}$. We remove all the rows and columns where the points in row $1\sim(i-1)$ are. In the rest of the figure, $|U'|=n-i+1$, $S'_3=\{\boldsymbol{a_1},\boldsymbol{a_p},\boldsymbol{a_q},\ldots,\boldsymbol{a_r}\}$ and the number of ways is $\frac{(n-i)!}{(p-1)(q-1)\cdot\ldots\cdot(r-1)}$. So,

$$\begin{aligned}g(S_1)&=\sum_{i=2}^{n+2-j}\frac{(n-i)!\,(n-j)!}{(p-1)(q-1)\cdot\ldots\cdot(r-1)(n-i-j+2)!}\\&=\frac{(n-1)!}{(p-1)(q-1)\cdot\ldots\cdot(r-1)(j-1)}\end{aligned}$$

So, for $\boldsymbol{C}=(0,1)$, change $S_3$ to $S_1$ by add $\boldsymbol{a_j}$ and we can get $g(S_1)=\frac{1}{j-1}g(S_3)$, where $j=2,3,\ldots,n$.

Then, for $k>1$, we discuss the relationship between $g(S_1)$ and $g(S_3)$. In $\boldsymbol{a_1}\sim\boldsymbol{a_n}$, let dimension *1*, $\{1,2,\ldots n\}$ and dimension *w+1*, $\{a_{1w},a_{2w},\ldots,a_{nw}\}$ form *k* groups of *n* *2*-dimensional vectors $(1,a_{1w}),(2,a_{2w}),\ldots(n,a_{nw})$ and $U^w=\{(1,a_{1w}),(2,a_{2w}),\ldots(n,a_{nw})\}$, where $w=1,2,\ldots,k$. Let $S^{w}\,\boldsymbol{OR},(0,1)\subseteq U^{w}$ and we can get the relationship between $S^w$ and $g(S^w)$ is the same as the relationship between $S$ and $g(S)$ when $\boldsymbol{C}=(0,1)$.

We define $S_{3,0}=\{\boldsymbol{a_1},\boldsymbol{a_p},\boldsymbol{a_q},\ldots,\boldsymbol{a_r}\}$ and there are 0 $(j,a_{jw})$ in all $S^w_{3,0}$. In *k+1*-dimensional vector, to change $S_{3,0}$ to $S_1$ by add $\boldsymbol{a_j}$, we need to change $l$ $S^w_{3,0}$ to $S^w_1$ by add $(j,a_{jw})$ if $c_l=1$. If $c_l=0$, we change $l$ $S^w_{3,0}$ to $S^w_1$ by add $(j,a_{jw})$ to change $S_{3,0}$ to $S_3$. We can get

$$g(S_1)=\boldsymbol{B}\cdot\boldsymbol{C}^T g(S_{3,0})$$
$$g(S_3)=\boldsymbol{B}\cdot\boldsymbol{C}'^T g(S_{3,0})$$

and

$$g(S_1)=\frac{\boldsymbol{B}\cdot\boldsymbol{C}^T}{\boldsymbol{B}\cdot\boldsymbol{C}'^T}g(S_3)=\frac{F_j(\boldsymbol{C})}{F_j(\boldsymbol{C}')}g(S_3)$$

We also get $\boldsymbol{a_1}$ in all $S^w$. So,

(a). If $c_k=1$,

$$O_{\boldsymbol{C}}(n,m)=\begin{cases}(n-1)!^k\sum\prod_{i=1}^{n-1}F_{j'_i}(\boldsymbol{C}'), & m=1\\ (n-1)!^k\sum\left(\prod_{i=1}^{m-1}F_{j_i}(\boldsymbol{C})\prod_{i=1}^{n-m}F_{j'_i}(\boldsymbol{C}')\right), & 2\le m\le n-1\\ (n-1)!^k\sum\prod_{i=1}^{n-1}F_{j_i}(\boldsymbol{C}), & m=n\\ 0, & m=0\text{ or }m>n\end{cases}$$

where $j_1,j_2,\cdots,j_{m-1}$ are all combinations consisting of $m-1$ elements in set $\{2,3,\ldots,n\}$ and $\{j'_1,j'_2,\cdots,j'_{n-m}\}=\{2,3,\ldots,n\}-\{j_1,j_2,\cdots,j_{m-1}\}$.

(b). If $c_k=0$,

$$O_{\boldsymbol{C}}(n,m)=\begin{cases}(n-1)!^k\sum\prod_{i=1}^{n-1}F_{j'_i}(\boldsymbol{C}'), & m=0\\ (n-1)!^k\sum(\prod_{i=1}^{m}F_{j_i}(\boldsymbol{C})\prod_{i=1}^{n-m-1}F_{j'_i}(\boldsymbol{C}')), & 1\le m\le n-2\\ (n-1)!^k\sum\prod_{i=1}^{n-1}F_{j_i}(\boldsymbol{C}), & m=n-1\\ 0, & m=-1 \text{ or } m>n-1\end{cases}$$

where $j_1,j_2,\cdots,j_m$ are all combinations consisting of $m$ elements in set $\{2,3,\dots,n\}$ and $\{j'_1,j'_2,\cdots,j'_{n-m-1}\}=\{2,3,\dots,n\}-\{j_1,j_2,\cdots,j_m\}$.

To sum up, for all $k,n\in N^+$ and $m\in N$,

$$O_{\boldsymbol{C}}(n,m+c_k-1)=\begin{cases}(n-1)!^k\sum(\prod_{i=1}^{m-1}F_{j_i}(\boldsymbol{C})\prod_{i=1}^{n-m}F_{j'_i}(\boldsymbol{C}')), & 1\le m\le n\\ 0, & otherwise\end{cases}$$

where $j_1,j_2,\cdots,j_{m-1}$ are all combinations consisting of $m-1$ elements in set $\{2,3,\dots,n\}$ and $\{j'_1,j'_2,\cdots,j'_{n-m}\}=\{2,3,\dots,n\}-\{j_1,j_2,\cdots,j_{m-1}\}$.

**Theorem 2.2.** *For all $k,n\in N^+$ and $m\ge c_k-1$, the recurrence formula of $\boldsymbol{C}$ sequential optimization numbers is*

$$O_{\boldsymbol{C}}(n+1,m+1)=n^kF_{n+1}(\boldsymbol{C})O_{\boldsymbol{C}}(n,m)+n^kF_{n+1}(\boldsymbol{C}')O_{\boldsymbol{C}}(n,m+1)$$

*and the boundary condition is*

$$O_{\boldsymbol{c}}(n,m)=\begin{cases}0, & m=c_k-1 \text{ or } m>c_k-1+n\\ 1, & m=c_k, n=1\end{cases}$$

**Proof of Theorem 2.2.** In Definition 2.2, for $c_k=1$, we divide the number of ways that $O_k(n+1,m+1)$ denotes into two parts, which are with and without $\boldsymbol{a_{n+1}}$, where $n>2$ and $m=2,3,\dots,n-1$.

First, ways with $\boldsymbol{a_{n+1}}$ and $m$ vectors from $\{\boldsymbol{a_1},\boldsymbol{a_2},\dots,\boldsymbol{a_n}\}$. let $j_m$ be $n+1$, $j_1,j_2,\cdots,j_{m-1}$ be all combinations consisting of *m-1* elements in set $\{2,3,\dots,n\}$ and $\{j'_1,j'_2,\cdots,j'_{n-m}\}=\{2,3,\dots,n\}-\{j_1,j_2,\cdots,j_{m-1}\}$.

$$\begin{aligned}n!^k\sum(\prod_{i=1}^{m}F_{j_i}(\boldsymbol{C})\prod_{i=1}^{n-m}F_{j'_i}(\boldsymbol{C}'))&=n!^k\,F_{n+1}(\boldsymbol{C})\sum(\prod_{i=1}^{m-1}F_{j_i}(\boldsymbol{C})\prod_{i=1}^{n-m}F_{j'_i}(\boldsymbol{C}'))\\&=n^kF_{n+1}(\boldsymbol{C})(n-1)!^k\sum(\prod_{i=1}^{m-1}F_{j_i}(\boldsymbol{C})\prod_{i=1}^{n-m}F_{j'_i}(\boldsymbol{C}'))\\&=n^kF_{n+1}(\boldsymbol{C})O_{\boldsymbol{C}}(n,m)\end{aligned}$$

Then, ways without $\boldsymbol{a_{n+1}}$ and m vectors from $\{\boldsymbol{a_1},\boldsymbol{a_2},\dots,\boldsymbol{a_n}\}$. let $j'_{n-m}$ be $n+1$, $j_1,j_2,\cdots,j_m$ are all combinations consisting of $m$ elements in set $\{2,3,\dots,n\}$ and $\{j'_1,j'_2,\cdots,j'_{n-m-1}\}=\{2,3,\dots,n\}-\{j_1,j_2,\cdots,j_m\}$.

$$n!^k\sum(\prod_{i=1}^{m}F_{j_i}(\boldsymbol{C})\prod_{i=1}^{n-m}F_{j'_i}(\boldsymbol{C}'))=n!^k\,F_{n+1}(\boldsymbol{C}')\sum(\prod_{i=1}^{m}F_{j_i}(\boldsymbol{C})\prod_{i=1}^{n-m-1}F_{j'_i}(\boldsymbol{C}'))$$

$$= n^k F_{n+1}(\boldsymbol{C}')(n-1)!^k \sum \left( \prod_{i=1}^{m} F_{j_i}(\boldsymbol{C}) \prod_{i=1}^{n-m-1} F_{j'_i}(\boldsymbol{C}') \right)$$

$$= n^k F_{n+1}(\boldsymbol{C}') O_{\boldsymbol{C}}(n, m+1)$$

To sum up,

$$O_{\boldsymbol{C}}(n+1, m+1) = n^k F_{n+1}(\boldsymbol{C}) O_{\boldsymbol{C}}(n, m) + n^k F_{n+1}(\boldsymbol{C}') O_{\boldsymbol{C}}(n, m+1)$$

We can get the boundary condition according to Theorem 2.1.

$$O_{\boldsymbol{C}}(n, m) = \begin{cases} 0, & m = c_k - 1 \text{ or } m > c_k - 1 + n \\ 1, & m = c_k, n = 1 \end{cases}$$

We can do the same thing for $c_k = 0$. We can prove that the recurrence formula works for all $k, n \in N^+$ and $m \geq c_k - 1$.

We give some properties and applications of **C** sequential optimization numbers below and we prove part of them in this paper.

**Lemma 2.1** *For* $c_k - 1 \leq m \leq n + c_k$,

$$O_{\boldsymbol{C}}(n, m) = O_{\boldsymbol{C}'}(n, n-m)$$

*For* $0 \leq m \leq n+1$,

$$O_{\boldsymbol{C}}(n, m + c_k - 1) = 0 = O_{\boldsymbol{C}'}(n, n + 1 - m + c'_k - 1).$$

**Proof of Lemma 2.1.** For $1 \leq t \leq n$,

$$O_{\boldsymbol{C}}(n, t + c_k - 1) = (n-1)!^k \sum \left( \prod_{i=1}^{t-1} F_{j_i}(\boldsymbol{C}) \prod_{i=1}^{n-t} F_{j'_i}(\boldsymbol{C}') \right)$$

$$= (n-1)!^k \sum \left( \prod_{i=1}^{n-t+1-1} F_{j'_i}(\boldsymbol{C}') \prod_{i=1}^{n-(n-t+1)} F_{j_i}(\boldsymbol{C}) \right)$$

$$= O_{\boldsymbol{C}'}(n, n - t + 1 + c'_k - 1)$$

$$= O_{\boldsymbol{C}'}(n, n - t + c'_k)$$

For $t = 0$ and $t = n + 1$,

$$O_{\boldsymbol{C}}(n, t + c_k - 1) = 0 = O_{\boldsymbol{C}'}(n, n + 1 - t + c'_k - 1)$$

Let $m = t + c_k - 1$,

$$O_{\boldsymbol{C}}(n, m) = O_{\boldsymbol{C}'}(n, n - m)$$

where $c_k - 1 \leq m \leq n + c_k$.

**Lemma 2.2.**

$$\sum_{m=0}^{n} O_{\boldsymbol{C}}(n, m) = n!^k$$

**Lemma 2.3.** *For* $\boldsymbol{C} = (0,1)$, $O_{\boldsymbol{C}}(n, m) = s_u(n, m)$, *where* $s_u(n, m)$ *are the unsigned Stirling numbers of first kind.*

**Lemma 2.4.** *Let* $\boldsymbol{C} = (0,1,1,\ldots,1)$ *be a k+1-tuple vector,* $O_{\boldsymbol{C}}(n, m) = O_k(n, m)$, *where* $O_k(n, m)$ *is k-dimensional sequential optimization numbers [Hui22].*

**Lemma 2.5.** *For all* $k, n \in N^+$, *we define* $x_{\boldsymbol{C}}^{1\uparrow} = x$,

$$x_{\boldsymbol{C}}^{n\uparrow} = x[F_2(\boldsymbol{C})x + F_2(\boldsymbol{C}')][2^k F_3(\boldsymbol{C})x + 2^k F_3(\boldsymbol{C}')] \cdots [(n-1)^k F_n(\boldsymbol{C})x + (n-1)^k F_n(\boldsymbol{C}')]$$

*where* $n \geq 2$ *and* $O_{\boldsymbol{C}}^{u}(n, m) = O_{\boldsymbol{C}}(n, m + c_k - 1)$. *Then, we can get*

$$x_{\boldsymbol{C}}^{n\uparrow} = \sum_{m=0}^{n} O_{\boldsymbol{C}}^{u}(n,m)\, x^m$$

*and the zero points are* $x = 0$ *and* $x = -\frac{F_m(\boldsymbol{C}')}{F_m(\boldsymbol{C})}$, *where* $m = 2,3,,\dots,n$. *We call* $O_{\boldsymbol{C}}^{u}(n,m)$ *unsigned* $\boldsymbol{C}$ *sequential optimization numbers.*

*For all* $k, n \in N^+$, *we define* $x_{\boldsymbol{C}}^{1\downarrow} = x$,

$$x_{\boldsymbol{C}}^{n\downarrow} = x[F_2(\boldsymbol{C})x - F_2(\boldsymbol{C}')][2^k F_3(\boldsymbol{C})x - 2^k F_3(\boldsymbol{C}')] \cdots [(n-1)^k F_n(\boldsymbol{C})x - (n-1)^k F_n(\boldsymbol{C}')]$$

*where* $n \geq 2$ *and* $O_{\boldsymbol{C}}^{s}(n,m) = (-1)^{n+m} O_{\boldsymbol{C}}(n, m + c_k - 1)$. *Then, we can get*

$$x_{\boldsymbol{C}}^{n\downarrow} = \sum_{m=0}^{n} O_{\boldsymbol{C}}^{s}(n,m)\, x^m$$

*and the zero points are* $x = 0$ *and* $x = \frac{F_m(\boldsymbol{C}')}{F_m(\boldsymbol{C})}$, *where* $m = 2,3,,\dots,n$. *We call* $O_{\boldsymbol{C}}^{s}(n,m)$ *signed* $\boldsymbol{C}$ *sequential optimization numbers.*

**Instance 2.1.** $\boldsymbol{C}$=$(c_0, c_1, c_2, \dots, c_k)$. *There are n boards and the sequence of them are fixed. Each board is painted one color and divided into k smaller boards. For all* $i = 1,2,\dots,k$, *the* $i$*th smaller boards of each board form a group and their height are* $1,2,\dots,n$ *respectively. In* $\{c_0, c_1, c_2, \dots, c_k\}$, $c_{p_1}, c_{p_2}, \dots, c_{p_q}$ *are all elements which are equal to 1. Following the direction of the arrow, the number of ways that m colors can be seen* $p_1$ *or* $p_2$ *or…or* $p_q$ *times is* $O_{\boldsymbol{C}}(n,m)$.

We call it $k$-dimensional color boards problem. An example is shown in Figure 2.a and the number of colors that can be seen is shown in Figure 2.b.

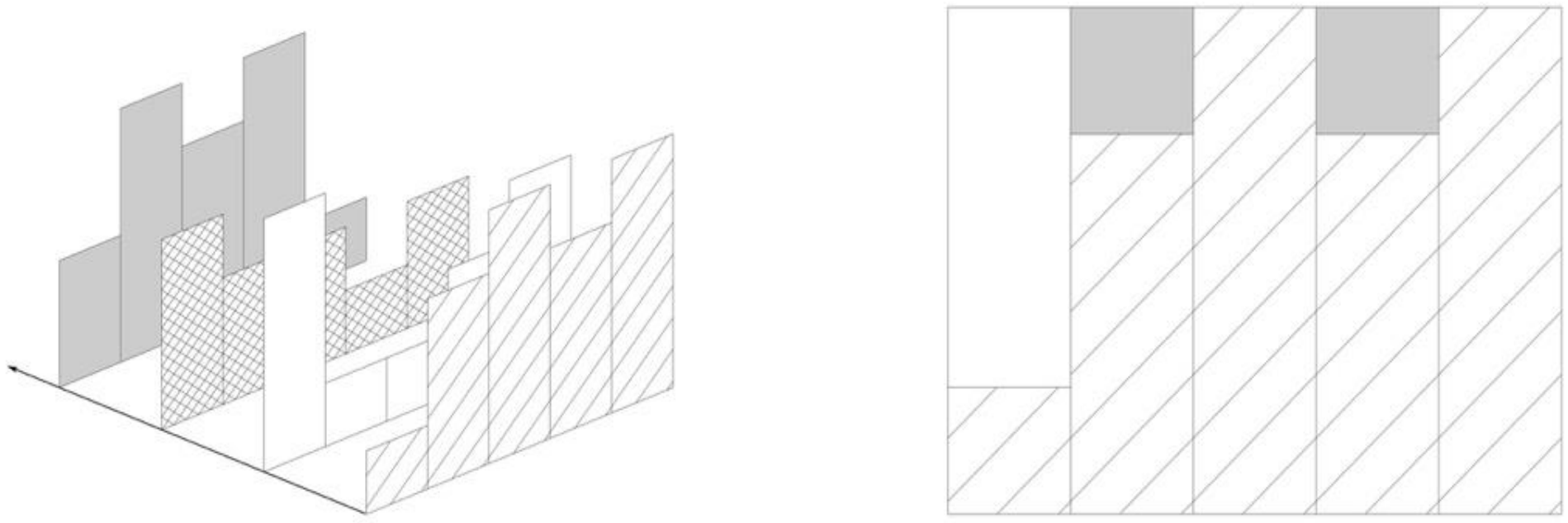

($a$) ($b$)

Figure 2: An example of $k$-dimensional color boards problem.

**Proof of Instance 2.1.** Label as the $i$th board follow direction of the arrow. Label the height of the smaller boards in same color from left to right as $h_{i1}, h_{i2}, \dots h_{ik}$ . Let $\boldsymbol{R} = (>, >, \dots, >)$ be $k$-dimensional relation vector, $\boldsymbol{h_i} = (i, h_{i1}, h_{i2}, \dots h_{ik})$ and $a_{ij} = n + 1 - h_{ij}$ , where $i = 1,2,\dots,n$ and $j = 1,2,\dots,k$. When $\boldsymbol{C} = (0,1)$, as shown in Figure 3, we can get same optimization set in Figure 1.$c$ follow direction of the line of sight. When $k \geq 1$, $j$th color is seen $l$ times means that change $l$ $S_{3,0}^{w}$ to $S_1^{w}$ by add $(j, a_{jw})$ in proof of Theorem 2.1. First color is seen $k$ time and $\boldsymbol{a_1}$ in all $S^w$ in proof of Theorem 2.1. So, we can conclude that the number of ways that $m$ colors can be seen $p_1$ or $p_2$ or…or $p_q$ times is $O_{\boldsymbol{C}}(n,m)$.

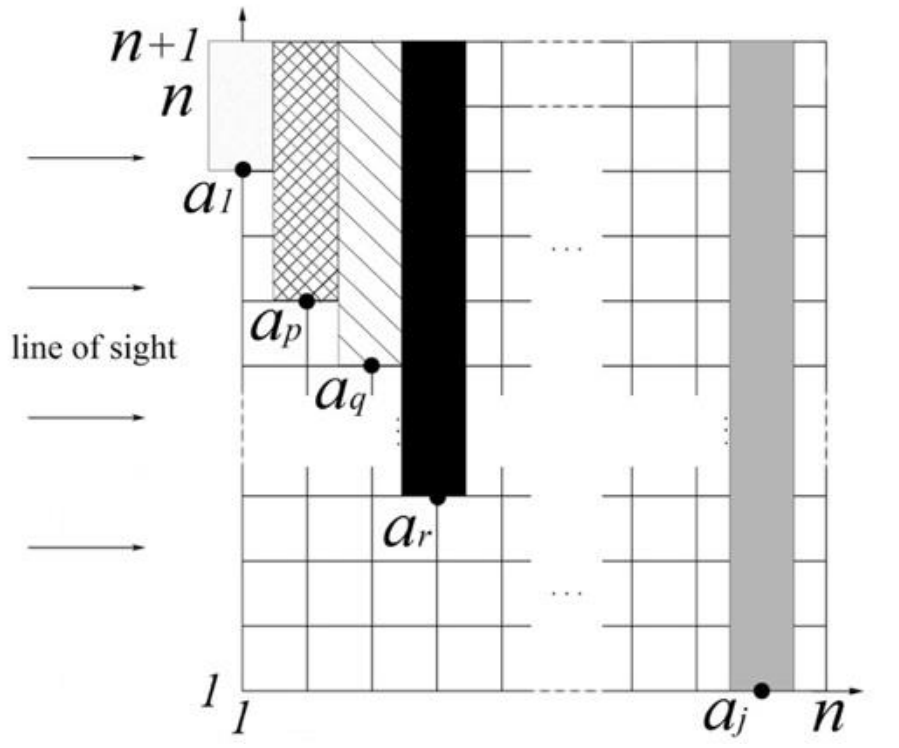

Figure 3: Relationship between smaller color boards group and optimization set.

The properties of **C** dimensional sequential optimization numbers above are basic properties of Stirling numbers of first kind when **C** is (0,1). Since Stirling numbers of first kind have been widely studied and a large number of achievements have been made, many of them can be transformed into the properties of **C** dimensional sequential optimization numbers and many valuable results can be obtained.

The properties of **C** dimensional sequential optimization numbers below are extended by properties of $k$-dimensional sequential optimization numbers [Hui22]. Theorem 2.3 is an upper bounder of **C** dimensional sequential optimization numbers. Theorem 2.4 and Lemma 2.6 show **C** dimensional sequential optimization numbers are almost concentrated in particular interval. Theorem 2.5 shows the upper ratio of upper bounder to themselves.

**Theorem 2.3.** *For all $k, n \in N^+$ and $m \in N$, we define*

$$O_{\boldsymbol{C}max}(n, m + c_k - 1) = \begin{cases} \dfrac{(n-1)!^k}{(m-1)!} (\boldsymbol{H}_n \cdot \boldsymbol{C}^T)^{m-1} \prod_{i=1}^{n-m} F_{i+1}(\boldsymbol{C}'), & 1 \le m \le n \\ 0, & otherwise \end{cases}$$

*where* $\boldsymbol{H}_n = (h_0, h_1, \dots, h_k)$ *and* $h_p = C_k^p \sum_{j=1}^{n-1} \frac{1}{j^p} = \sum_{j=2}^{n} C_k^p \frac{1}{(j-1)^p} = \sum_{j=2}^{n} b_{j,p}$ . *We can get* $O_{\boldsymbol{C}max}(n, m + c_k - 1) \ge O_{\boldsymbol{C}}(n, m + c_k - 1)$.

We can get a formula below.

$$\boldsymbol{H}_n \cdot \boldsymbol{C}^T = \sum_{p=1}^{k} h_p c_p = \sum_{p=1}^{k} \sum_{j=2}^{n} b_{j,p}\, c_p = \sum_{j=2}^{n} \sum_{p=1}^{k} b_{j,p}\, c_p = \sum_{j=2}^{n} F_j(\boldsymbol{C}) \tag{4}$$

**Proof of Theorem 2.3.** We prove it using Theorem 2.2. First, we prove the boundary condition. For all $k, n \in N^+$,

$$O_{\boldsymbol{C}max}(n, \mathrm{m} + c_k - 1) = \begin{cases} (n-1)!^k \prod_{i=1}^{n-1} F_{i+1}(\boldsymbol{C}'), & m = 1 \\ 0, & m = 0 \text{ or } m > n \end{cases}$$

So, the boundary condition satisfy $O_{\boldsymbol{C}max}(n, m + c_k - 1) \ge O_{\boldsymbol{C}}(n, m + c_k - 1)$.

Then, we prove the inequality for $2 \le m \le n$ by mathematical induction.

Basis step: For $n = 2$, according to formula (2),

$$O_{\boldsymbol{C}max}(2,2 + c_k - 1) = \boldsymbol{H}_2 \cdot \boldsymbol{C}^T = F_2(\boldsymbol{C}) = O_{\boldsymbol{C}}(2,2 + c_k - 1) \qquad formula\ (4)$$

So, for $n = 2$, $O_{\boldsymbol{C}max}(n, m + c_k - 1) \geq O_{\boldsymbol{C}}(n, m + c_k - 1)$.

Inductive step: For $n \geq 2$ and $2 \leq m \leq n$, we assume $O_{\boldsymbol{C}max}(n, m + c_k - 1) \geq O_{\boldsymbol{C}}(n, m + c_k - 1)$. Then, we prove $O_{\boldsymbol{C}max}(n + 1, m + c_k - 1) \geq O_{\boldsymbol{C}}(n + 1, m + c_k - 1)$.

$$O_{\boldsymbol{C}max}(n + 1, \mathrm{m} + c_k - 1) = \frac{n!^k}{(m-1)!}(\boldsymbol{H}_{n+1} \cdot \boldsymbol{C}^T)^{m-1} \prod_{i=1}^{n+1-\mathrm{m}} F_{i+1}(\boldsymbol{C}')$$

$$= \frac{n!^k}{(m-1)!}[\boldsymbol{H}_n \cdot \boldsymbol{C}^T + F_{n+1}(\boldsymbol{C})]^{m-1} \prod_{i=1}^{n+1-\mathrm{m}} F_{i+1}(\boldsymbol{C}') \qquad formula\ (4)$$

$$\geq \frac{n!^k}{(m-1)!} C_{m-1}^0 (\boldsymbol{H}_n \cdot \boldsymbol{C}^T)^{m-1} \prod_{i=1}^{n+1-\mathrm{m}} F_{i+1}(\boldsymbol{C}') + \frac{n!^k}{(m-1)!} C_{m-1}^1 F_{n+1}(\boldsymbol{C})(\boldsymbol{H}_n \cdot \boldsymbol{C}^T)^{m-2} \prod_{i=1}^{n+1-\mathrm{m}} F_{i+1}(\boldsymbol{C}')$$

$$= \frac{n^k (n-1)!^k}{(m-1)!}(\boldsymbol{H}_n \cdot \boldsymbol{C}^T)^{m-1} \prod_{i=1}^{n+1-\mathrm{m}} F_{i+1}(\boldsymbol{C}') + \frac{n^k (n-1)!^k}{(m-2)!} F_{n+1}(\boldsymbol{C})(\boldsymbol{H}_n \cdot \boldsymbol{C}^T)^{m-2} \prod_{i=1}^{n+1-\mathrm{m}} F_{i+1}(\boldsymbol{C}')$$

$$= F_{n+2-\mathrm{m}}(\boldsymbol{C}') \frac{n^k (n-1)!^k}{(m-1)!}(\boldsymbol{H}_n \cdot \boldsymbol{C}^T)^{m-1} \prod_{i=1}^{n-\mathrm{m}} F_{i+1}(\boldsymbol{C}') + \frac{n^k (n-1)!^k}{(m-2)!} F_{n+1}(\boldsymbol{C})(\boldsymbol{H}_n \cdot \boldsymbol{C}^T)^{m-2} \prod_{i=1}^{n+1-\mathrm{m}} F_{i+1}(\boldsymbol{C}')$$

$$\geq F_{n+1}(\boldsymbol{C}') \frac{n^k (n-1)!^k}{(m-1)!}(\boldsymbol{H}_n \cdot \boldsymbol{C}^T)^{m-1} \prod_{i=1}^{n-\mathrm{m}} F_{i+1}(\boldsymbol{C}') + \frac{n^k (n-1)!^k}{(m-2)!} F_{n+1}(\boldsymbol{C})(\boldsymbol{H}_n \cdot \boldsymbol{C}^T)^{m-2} \prod_{i=1}^{n+1-\mathrm{m}} F_{i+1}(\boldsymbol{C}') \qquad formula\ (1)$$

$$= n^k F_{n+1}(\boldsymbol{C}') O_{\boldsymbol{C}max}(n, m + c_k - 1) + n^k F_{n+1}(\boldsymbol{C}) O_{\boldsymbol{C}max}(n, m - 1 + c_k - 1)$$

$$\geq n^k F_{n+1}(\boldsymbol{C}') O_{\boldsymbol{C}}(n, m + c_k - 1) + n^k F_{n+1}(\boldsymbol{C}) O_{\boldsymbol{C}}(n, m - 1 + c_k - 1)$$

$$= O_{\boldsymbol{C}}(n + 1, \mathrm{m} + c_k - 1)$$

We can do same thing for $m = n + 1$. To sum up, for all $k, n \in N^+$ and $m \in \mathrm{N}$, $O_{\boldsymbol{C}max}(n, m + c_k - 1) \geq O_{\boldsymbol{C}}(n, m + c_k - 1)$.

**Theorem 2.4.** *For all* $n \geq 2$, $k \geq 1$ *and* $c_0 = 0$, $Pr\ [O_{\boldsymbol{C}}(n, m + c_k - 1)] = \frac{O_{\boldsymbol{C}}(n, m + c_k - 1)}{n!^k}$ *be the probability of* $O_{\boldsymbol{C}}(n, m + c_k - 1)$ *and* $M = \lceil ekc_1[log(n-1) + 1] + e\frac{\pi^2}{6}\sum_{p=2}^{k} c_p C_k^p \rceil + M_1$, *where* $M_1$ *is a positive integer. We can get*

$$Pr\ [O_{\boldsymbol{C}}(n, m > M + c_k - 1)] \leq e^{-M_1}$$

**Proof of Theorem 2.4.** Let $Pr\ [O_{\boldsymbol{C}max}(n, m + c_k - 1)] = \frac{O_{\boldsymbol{C}max}(n, m + c_k - 1)}{n!^k}$ be the probability of the upper bound of **C** sequential optimization numbers. For $c_0 = 0$ and $1 \leq m \leq n - 1$,

$$Pr\ [O_{kmax}(n, m + c_k)] = \frac{O_{kmax}(n, m + c_k)}{n^k}$$

$$= \frac{1}{n^k m!} (\boldsymbol{H}_n \cdot \boldsymbol{C}^T)^m \prod_{i=1}^{n-m-1} F_{i+1}(\boldsymbol{C}')$$

$$\leq \frac{1}{n^k m!} (\boldsymbol{H}_n \cdot \boldsymbol{C}^T)^m (n - m)^k \qquad formula\ (2)$$

$$\leq \frac{1}{m!} (\boldsymbol{H}_n \cdot \boldsymbol{C}^T)^m$$

and

$$\frac{Pr\ [O_{\boldsymbol{C}max}(n, m + c_k)]}{Pr\ [O_{\boldsymbol{C}max}(n, m + c_k - 1)]} = \frac{\boldsymbol{H}_n \cdot \boldsymbol{C}^T}{m F_{n-m+1}(\boldsymbol{C}')}$$

$$\leq \frac{\boldsymbol{H}_n \cdot \boldsymbol{C}^T}{m} \qquad formula\ (3)$$

Since $h_p = \sum_{j=2}^{n} C_k^p \frac{1}{(j-1)^p}$, we can get $h_0 = n - 1$, $h_1 = k \sum_{j=2}^{n} \frac{1}{j-1} \leq k[log(n - 1) + 1]$ and for $p \geq 2$, $h_p \leq C_k^p \sum_{j=2}^{n} \frac{1}{(j-1)^2} \leq \frac{\pi^2}{6} C_k^p$. We know $m! > \sqrt{2\pi m}(\frac{m}{e})^m$ (Stirling's approximation).

$$\frac{Pr\ [O_{\boldsymbol{C}max}(n, m + c_k)]}{Pr\ [O_{\boldsymbol{C}max}(n, m + c_k - 1)]} \leq \frac{c_1\, C_k^1[\log(n - 1) + 1] + \frac{\pi^2}{6} \sum_{p=2}^{k} c_p C_k^p}{m}$$

When $m \geq \lceil \mathrm{ek}c_1[log(n - 1) + 1] + e\frac{\pi^2}{6} \sum_{p=2}^{k} c_p C_k^p \rceil$,

$$\frac{Pr\ [O_{\boldsymbol{C}max}(n, m + c_k)]}{Pr\ [O_{\boldsymbol{C}max}(n, m + c_k - 1)]} \leq \frac{1}{e}$$

and

$$Pr\ [O_{\boldsymbol{C}max}(n, m + c_k)] = \frac{1}{m!} (\boldsymbol{H}_n \cdot \boldsymbol{C}^T)^m$$

$$\leq \frac{\{c_1\, C_k^1[\log(n - 1) + 1] + \frac{\pi^2}{6} \sum_{p=2}^{k} c_p C_k^p\}^m}{\sqrt{2\pi m}(\frac{m}{e})^m}$$

$$= \frac{1}{\sqrt{2\pi m}} [\frac{\mathrm{ek}c_1[log(n - 1) + 1] + e\frac{\pi^2}{6} \sum_{p=2}^{k} c_p C_k^p}{m}]^m$$

$$\leq e^{-1}$$

When $M = \lceil \mathrm{ek}c_1[log(n - 1) + 1] + e\frac{\pi^2}{6} \sum_{p=2}^{k} c_p C_k^p \rceil + M_1$, where $M_1$ is a positive integer.

$$Pr\ [O_{\boldsymbol{C}}(n, M + c_k - 1)] \leq e^{-M_1}$$

So,

$$Pr\,[O_{\boldsymbol{C}}(n, m > M + c_k - 1)] \leq Pr\,[O_{\boldsymbol{C}max}(n, m > M + c_k - 1)]$$

$$= \sum_{i=M+1}^{n} Pr\,[O_{\boldsymbol{C}max}(n, i)]$$

$$\leq e^{-M_1}$$

**Lemma 2.6** *If* $c_0 = 1$, *we can get* $c'_0 = 0$ *and* $P\,r[O_{\boldsymbol{C}}(n, m > n - M + c_k)] = P\,r[O_{\boldsymbol{C}'}(n, m < M + c'_k - 1)] \leq e^{-M_1}$, *where* $M = \lceil ekc'_1[log(n-1)+1] + e\frac{\pi^2}{6}\sum_{p=2}^{k} c'_p C_k^p \rceil + M_1$ *and* $M_1$ *is a positive integer.*

**Proof of Lemma 2.6.** $c'_0 = 1 - c_0 = 0$.

$$P\,r[O_{\boldsymbol{C}'}(n, m < M + c'_k - 1)] \leq e^{-M_1}$$

$$P\,r[O_{\boldsymbol{C}}(n, m > n + 1 - M + c_k - 1)] = P\,r[O_{\boldsymbol{C}'}(n, m < M + c'_k - 1)] \leq e^{-M_1} \qquad Lemma\ 2.1$$

where $M = \lceil ekc'_1[log(n-1)+1] + e\frac{\pi^2}{6}\sum_{p=2}^{k} c'_p C_k^p \rceil + M_1$ and $M_1$ is a positive integer.

**Theorem 2.5.** *For all* $k \geq 1$, $n \geq 2$ *and* $0 \leq m \leq n$, $O_{\boldsymbol{C}}(n, m) \leq O_{\boldsymbol{C}'max}(n, n-m)$,

$$\frac{\sum_{m=0}^{n} O_{\boldsymbol{C}max}(n, m)}{\sum_{m=0}^{n} O_{\boldsymbol{C}}(n, m)} \leq e^{\lambda}$$

*and*

$$\frac{\sum_{m=0}^{n} O_{\boldsymbol{C}'max}(n, m)}{\sum_{m=0}^{n} O_{\boldsymbol{C}}(n, m)} \leq e^{\lambda'}$$

*where* $\lambda = c_0(n-1) + c_1 k[log(n-1)+1] + \frac{\pi^2}{6}\sum_{p=2}^{k} c_p C_k^p$ *and* $\lambda' = c'_0(n-1) + c'_1 k[log(n-1)+1] + \frac{\pi^2}{6}\sum_{p=2}^{k} c'_p C_k^p$.*In particular, for* $s_u(n, m) = O_{(0,1)}(n, m)$,

$$\frac{\sum_{m=0}^{n} s_{umax}(n, m)}{\sum_{m=0}^{n} s_u(n, m)} \leq \frac{(n-1)}{n} e^{\sum_{j=2}^{n}\frac{1}{j-1} - log\,(n-1)} \leq e^{\gamma}$$

*where* $\gamma$ *is Euler-Mascheroni constant and* $e^{\gamma} \leq 1.7811$.

**Proof of Theorem 2.5.**

$$O_{\boldsymbol{C}}(n, m) = O_{\boldsymbol{C}'}(n, n-m) \leq O_{\boldsymbol{C}'max}(n, n-m) \qquad Lemma\ 2.2$$

Let

$$\lambda = \boldsymbol{H}_n \cdot \boldsymbol{C}^T$$

$$\leq c_0(n-1) + c_1 k[\log(n-1)+1] + \frac{\pi^2}{6}\sum_{p=2}^{k} c_p C_k^p$$

and

$$\lambda' = \boldsymbol{H}_n \cdot \boldsymbol{C}'^T$$

$$\leq c'_0(n-1) + c'_1 k[\log(n-1)+1] + \frac{\pi^2}{6}\sum_{p=2}^{k} c'_p C_k^p$$

then,

$$\frac{\sum_{m=0}^{n} O_{\boldsymbol{C}max}(n,m)}{\sum_{m=0}^{n} O_{\boldsymbol{C}}(n,i)} = \frac{1}{n!^k}\sum_{m=1}^{n}\frac{(n-1)!^k}{(m-1)!}(\boldsymbol{H}_n\cdot\boldsymbol{C}^T)^{m-1}\prod_{i=1}^{n-\mathrm{m}} F_{i+1}(\boldsymbol{C}')$$

$$\leq \frac{1}{n^k}\sum_{m=1}^{n}\frac{\lambda^{i-1}}{(m-1)!}(n-m+1)^k$$

$$\leq \sum_{m=1}^{n}\frac{\lambda^{i-1}}{(m-1)!}$$

$$\leq e^{\lambda} \qquad (Taylor\ theorem)$$

and

$$\frac{\sum_{m=0}^{n} O_{\boldsymbol{C}'max}(n,m)}{\sum_{m=0}^{n} O_{\boldsymbol{C}}(n,m)} = \frac{1}{n!^k}\sum_{m=1}^{n}\frac{(n-1)!^k}{(m-1)!}(\boldsymbol{H}_n\cdot\boldsymbol{C}'^T)^{m-1}\prod_{i=1}^{n-\mathrm{m}} F_{i+1}(\boldsymbol{C})$$

$$\leq \frac{1}{n^k}\sum_{m=1}^{n}\frac{\lambda'^{i-1}}{(m-1)!}(n-m+1)^k$$

$$\leq \sum_{m=1}^{n}\frac{\lambda'^{i-1}}{(m-1)!}$$

$$\leq e^{\lambda'} \qquad (Taylor\ theorem)$$

In particular, for $s(n,m) = O_{(0,1)}(n,m)$,

$$\frac{\sum_{m=1}^{n} s_{umax}(n,m)}{\sum_{m=1}^{n} s_u(n,m)} = \frac{1}{n!}\sum_{m=1}^{n}\frac{(n-1)!}{(m-1)!}(\boldsymbol{H}_n\cdot\boldsymbol{C}^T)^{m-1}\prod_{i=1}^{n-\mathrm{m}} F_{i+1}(\boldsymbol{C}')$$

$$\leq \frac{1}{n}\sum_{m=1}^{n}\frac{\lambda^{i-1}}{(m-1)!}$$

$$\leq \frac{e^{\lambda}}{n} \qquad (Taylor\ theorem)$$

$$= \frac{1}{n}e^{\sum_{j=2}^{n}\frac{1}{j-1}}$$

Let

$$a_n = \frac{1}{n}e^{\sum_{j=2}^{n}\frac{1}{j-1}}$$

and

$$f(x) = \frac{x}{x+1}e^{\frac{1}{x}}$$

where $n > 0$ and $x > 0$.

$$\frac{a_{n+1}}{a_n} = \frac{n}{n+1} e^{\frac{1}{n}}$$

$$f'(x) = -\frac{1}{x(x+1)^2} e^{\frac{1}{x}} < 0$$

$$\lim_{x\to+\infty} f(x) = 1$$

So, $f(x) > 1, \frac{a_{n+1}}{a_n} > 1$ and

$$\begin{aligned}\frac{\sum_{m=1}^{n} s_{umax}(n,m)}{\sum_{m=1}^{n} s_u(n,m)} &\leq \lim_{n\to+\infty} \frac{1}{n} e^{\sum_{j=2}^{n}\frac{1}{j-1}} \\ &= \lim_{n\to+\infty} \frac{(n-1)}{n} e^{\sum_{j=2}^{n}\frac{1}{j-1} - log\,(n-1)} \\ &= e^{\gamma}\end{aligned}$$

where $\gamma$ is Euler-Mascheroni constant and $e^{\gamma} \leq 1.7811$.

## 3 CONCLUSION

In the article [Hui22], the proof of the Stirling numbers of first kind is more basic, which led to the discovery of **C** sequential optimization numbers. We deal with many inequalities succinctly and roughly and some of them still have much room for improvement. For the *k+1*-tuple vector **C**, there are $2^k$ kinds of sequences and some of them can have special properties.